\documentclass[12pt,twoside,a4paper]{amsart}
\usepackage{amssymb}
\date{\today}

\def\End{{\rm End}}
\def\ann{{\rm ann}}
\def\deg{\text{deg}\,}

\def\depth{\text{depth}\,}

\def\w{\wedge}

\def\dbar{\bar\partial}

\def\R{{\mathbb R}}
\def\C{{\mathbb C}}
\def\w{{\wedge}}
\def\P{{\mathbb P}}

\def\D{{\mathcal D}}

\def\T{{\mathcal T}}
\def\S{{\mathcal S}}

\def\J{{\mathcal J}}
\def\F{{\mathcal F}}

\def\Hom{{\rm Hom\, }}
\def\codim{{\rm codim\,}}

\def\Im{{\rm Im\, }}

\def\Ker{{\rm Ker\,  }}
\def\rank{{\rm rank\, }}

\def\Z{{\mathbb Z}}
\def\E{{\mathcal E}}

\def\Ok{{\mathcal O}}

\def\Re{{\rm Re\,  }}

\def\U{{\mathcal U}}

\def\be{\begin{equation}}
\def\ee{\end{equation}}

\newtheorem{thm}{Theorem}[section]
\newtheorem{lma}[thm]{Lemma}
\newtheorem{cor}[thm]{Corollary}
\newtheorem{prop}[thm]{Proposition}

\theoremstyle{definition}

\theoremstyle{remark}

\newtheorem{preremark}{Remark}
\newtheorem{preex}{Example}

\newenvironment{remark}{\begin{preremark}}{\qed\end{preremark}}
\newenvironment{ex}{\begin{preex}}{\qed\end{preex}}

\numberwithin{equation}{section}

\title[Residue currents with prescribed  annihilator ideals]%%and dualityoetherian residue currents]%
{Residue currents with prescribed  annihilator ideals}%% and the fundamental prindiple}

\begin{document}

\date{\today}

\author{Mats Andersson \& Elizabeth Wulcan}

\address{Department of Mathematics\\Chalmers University of Technology and the University of 
G\"oteborg\\S-412 96 G\"OTEBORG\\SWEDEN}

\email{matsa@math.chalmers.se,   wulcan@math.chalmers.se}

\subjclass{32A26, 32A27, 32C35}

\thanks{The first author was
  partially supported by the Swedish Natural Science
  Research Council}

\begin{abstract}
Given a coherent  ideal sheaf $J$ 
we construct locally a vector-valued residue current $R$ whose annihilator is
precisely the given sheaf. In case $J$ is a complete intersection,
$R$ is just the classical Coleff-Herrera product. By means of these currents
we can extend various results, previously known for a complete intersection,
to general ideal sheaves. Combining  with 
 integral formulas  we obtain a residue version of the Ehrenpreis-Palamodov
fundamental principle. 
%%Analogous  results hold true also  for a coherent subsheaf of a locally free 
%%analytic sheaf.
\end{abstract}

%%

%%uppsats/integralformler2/genexact.....

%%%%%%%%%%

\maketitle

\section{Introduction}
%%%
Let $h=h_1,\ldots,h_m$ be a tuple of holomorphic functions 
such that their common zero set $Z$ has codimension $m$,
and let
\begin{equation}\label{ch}
\mu^h=\dbar\frac{1}{h_1}\w\ldots\w\dbar\frac{1}{h_m}
\end{equation}
be the  Coleff-Herrera product  introduced in  \cite{CH}.
Dickenstein-Sessa, \cite{DS}, and Passare, \cite{P1},
independently proved the duality principle,
that a holomorphic function $\phi$ is in the ideal
sheaf $\J(h)$ generated by $h_1,\ldots,h_m$ if and only if  
the current $\phi \mu^h$ vanishes, i.e., $\phi$
belongs to the annihilator $\ann \mu^h$.
%%%%%%%%%%%%%%%%%
Given  any coherent ideal sheaf $\J$ one can locally find  a finite tuple
$\gamma=(\gamma_1,\ldots\gamma_\mu)$  
of so-called Coleff-Herrera currents  such that 
$\J=\ann\gamma=\cap_j \ann \gamma_j$; 
this is closely related to the existence of Noetherian operators,
see \cite{JEB4}.  
However, much of the utility of  the duality principle depends on the
fact that the current $\mu^h$  fits into various  division-interpolation
integral formulas,
see,  e.g., \cite{BB}, \cite{P1}, \cite{BP}, \cite{BY0}, and \cite{BGVY}. 
Therefore it is natural to  look  for an analogue, for a general
ideal sheaf, with this extra property.

To begin with we consider  an  arbitrary
complex of Hermitian holomorphic vector bundles over a complex manifold $X$,
\begin{equation}\label{plex}
0\to E_N\stackrel{f_N}{\longrightarrow}\ldots\stackrel{f_3}{\longrightarrow} 
E_2\stackrel{f_2}{\longrightarrow}
E_1\stackrel{f_1}{\longrightarrow}E_0,
\end{equation}
that is exact outside an analytic  variety $Z$ of positive codimension.
To this complex $E_\bullet$ we associate  %%%in Section~\ref{hungrig} 
a current $R=R(E_\bullet)$ taking  values in  $\End(\oplus_k E_k)$
and with support on $Z$. This current in a certain way 
 measures the lack of exactness of the associated  complex of locally
free sheaves of $\O$-modules of sections of $E_k$
\begin{equation}\label{splex2}
0\to \Ok(E_N)\to \cdots\to \Ok(E_1)\to\Ok(E_0).
\end{equation}
%%%%
Let $R^\ell$ denote the component of $R$ that takes values in $\Hom(E_\ell,\oplus_k E_k)$.
It turns out that   \eqref{splex2} is exact if and only if $R^\ell=0$ for $\ell\ge 1$
(Theorem~\ref{main1}).  
Let $\J=\Im(\Ok(E_1)\to\Ok(E_0))$.
The  main result in this paper is the following:

\begin{thm}\label{skalle}
Suppose that the sheaf complex \eqref{splex2} is exact. Then the associated residue current
$R$ has its support on the set $Z$
where the sheaf  $\Ok(E_0)/\J$ is not locally free,
and a local holomorphic section $\phi$ of
$E_0$ is in $\J$ if and only if $\phi$ is generically in the image of $f_1$ and
the residue current $R\phi$ vanishes.
\end{thm}

The set $Z$ is precisely the set where the mapping $f_1$ does not have optimal rank.
If $f_1$ is generically surjective, or equivalently 
$\ann (\Ok(E_0)/\J)$ is nonzero, thus $\phi\in \J$ if and only if $R\phi=0$.
In  this case $Z$ is the zero locus of $\ann (\Ok(E_0)/\J)$.
In particular as soon as $\J$ is a nontrivial ideal sheaf ($\rank E_0=1$)
then $R$ has its support on the zero locus of $\J$ and $\phi\in \J$ if and only if
$R\phi=0$.
In analogy with Noetherian differential operators it is natural  to say that
$R$ is a Noetherian residue current for $\J$.

If $\J$ is any  coherent subsheaf of some  locally free sheaf $\Ok(E_0)$,
then at least locally   $\Ok(E_0)/\J$ admits a resolution \eqref{splex2},
and if we equip the corresponding complex of vector bundles with
any Hermitian metric we thus locally get a current $R$ 
as in Theorem~\ref{skalle}. In case $\J$ is defined by a complete intersection,
the Koszul complex provides a resolution, and the resulting  residue current is just the 
Coleff-Herrera product, see Example~\ref{ex1} below.
In general it is just as hard to find resolutions of ideals as to find, e.g., Noetherian
differential  operators, so Theorem~\ref{skalle} will not  contribute  to 
effectivity questions, but it turns out to be useful in several other ways.

\smallskip

If $\Ok(E_0)/\J$ is a sheaf of Cohen-Macaulay modules, 
the associated    current $R$ is independent of the Hermitian metrics and it is 
essentially canonical, see Section~\ref{cmsection}  for precise statements. 
%%%
In the Cohen-Macaulay  case we can also define a cohomological residue for $\J$,
so that the  cohomological duality principle for
a complete intersection ideal extends (Theorem~\ref{dualthm}).

\smallskip
%%%
Combined with the framework of integral formulas developed in \cite{A7}, we 
present in Section~\ref{intf}  a holomorphic  division  formula,
\eqref{decomp},  for sections of $E_k$. 
In particular, as soon as  $\phi\in\J$,
this formula  provides  an explicit realization of the membership. 
By   a similar  integral formula we  obtain a residue characterization 
(Theorem~\ref{smoothy}) of the sheaf   $\E \J$
of $\E$-modules generated by $\J$.

\smallskip
Given a module  $J$  over $\C[z_1,\ldots,z_n]$, 
generated  by an   $r_0\times r_1$-matrix $F(z)$ of polynomials in $\C^n$ of 
generic rank $r_0$
we can find a global Noetherian residue current $R$ for the corresponding
sheaf $\J$ in $\C^n$.
It  is obtained  from  a resolution  of the  module 
 over  the graded ring  $\C[z_0,\ldots, z_n]$  induced by a homogenization of $F$.
We can use this current to prove a generalization of  Max Noether's
classical $AF+BG$ theorem.
Our main application  is a  residue version of the general
fundamental principle:
If $F^T$ is the transpose of $F$, then any smooth
solution to  $F^T(i\partial/\partial t)\xi=0$ on a smoothly bounded convex set 
in $\R^n$ can be written 
$$
\xi(t)=\int_{\C^n}R^T(\zeta) A(\zeta) e^{-i\langle t,\zeta\rangle},
$$
for an  appropriate (explicitly given matrix of smooth functions) $A$;
here $R^T$ is the transpose of $R$. Conversely, 
since $R$ is  Noetherian,  any
$\xi(t)$ given in this way is a homogeneous solution.
This follows %%from an appropriately weighted integral formula, 
along the same lines as in \cite{BP}, where
this result was obtained  for a complete intersection
$F$  by means of the Coleff-Herrera product.

\smallskip

Throughout this paper, 
$\E_\bullet(E)$, $\D_\bullet(E)$, $\D'_\bullet(E)$, and
$\Ok(E)$ denote the sheaves  of smooth forms, test forms,
currents, and holomorphic functions, respectively with values in the
vector bundle $E$.

{\bf Acknowledgement:} We express our sincere gratitude to  
Jan-Erik Bj\"ork, Ralf  Fr\"oberg,  and Alain Yger for invaluable
discussions on these matters.  We also thank the referee for
several suggestions for improvements of the presentation.

\section{Residue currents of generically exact complexes}\label{hungrig}

Let $E,Q$ be Hermitian holomorphic vector bundles over a connected manifold $X$
and let $f\colon E\to Q$ be  a holomorphic  morphism. 
If $f$ has optimal rank $\rho$ then the rank is precisely $\rho$ outside the
analytic set  $Z=\{F=0\}$, where $F=\det^\rho f$ is a section of
$\Lambda^\rho E^*\otimes \Lambda^\rho Q$.
Let $\sigma\colon Q\to E$ be  the minimal inverse in $X\setminus Z$,  
 i.e., $\sigma\xi$ is the
minimal solution to $f\eta=\xi$ if $\xi$ is in the image of $f$
and $\sigma\xi=0$ if $\xi$ is orthogonal to $\Im f$.
%%%
Then clearly $\sigma$ is smooth outside $Z$,   and following
the proof of  Lemma~4.1 in \cite{A5}  we get
%%If $s$ be the section of $E\otimes Q^*$ that is dual to  $f$, 
%%then  $S=\det_q s$ is the  section of
%%$\Lambda^qE\otimes\Lambda^qQ^*$ that is dual to $F$, and 
%%\begin{equation}\label{cl2}
%%\sigma=(\delta_f)_{q-1}S/|F|^2.
%%\end{equation}
%%Clearly $\sigma$ is smooth outside $Z$. We also have

\begin{lma}\label{sing}
If $F=F^0F'$ in $X$, where $F^0$ is a holomorphic function  and
$F'$ is non-vanishing, then
$F^0\sigma$  is smooth across $Z$.
\end{lma}

Let 
\begin{equation}\label{complex}
0\to E_N\stackrel{f_N}{\longrightarrow}E_{N-1}\stackrel{f_{N-1}}{\longrightarrow}
\ldots \stackrel{f_{-M+2}}{\longrightarrow}
E_{-M+1}\stackrel{f_{-M+1}}{\longrightarrow}E_{-M}\to 0
\end{equation}
be a holomorphic complex of Hermitian vector bundles 
over the  $n$-dimensional complex manifold $X$,
and assume that it is 
pointwise exact outside the  analytic set $Z$ of positive codimension.
Then  for each $k$, 
$\rank f_k$  is constant  in  $X\setminus Z$  and  equal to
\begin{equation}\label{rhodef}
\rho_k=\dim E_{k}-\dim E_{k+1}+\cdots \pm \dim E_N.
\end{equation}
The bundle $E=\oplus E_k$ has a natural superbundle structure, i.e., a $\Z_2$-grading,
$E=E^+\oplus E^-$, $E^+$ and $E^-$ being the subspaces of even and odd elements,
respectively,  by letting
$E^+=\oplus _{2k}E_k$ and $E^-=\oplus_{2k+1}E_k$, see \cite{Qu} and, e.g.,  
\cite{A7},  for details.
The mappings $f=\sum f_j$  and $\dbar$ are then odd mappings on 
$\D'_\bullet(E)$  %%%\D'_\bullet(X)\otimes_{\E(X)}\E(X,E),$
and  they anticommute so that $\nabla^2=0$, where $\nabla=f-\dbar$ is 
(minus) the $(0,1)$-part of Quillen's superconnection $D-\dbar$.
Moreover, $\nabla$  extends to an odd mapping  $\nabla_\End$ on $\D'_\bullet(\End E)$ 
and $\nabla_{\End}^2=0$.
%%%%
%%
In  $X\setminus Z$ let $\sigma_k\colon E_{k-1}\to E_k$ 
be the minimal inverses of $f_k$. 
If $\sigma=\sigma_{-M+1}+\cdots+\sigma_N\colon E\to E$ 
and  $I$ denotes  the identity endomorphism on $E$, then 
$f\sigma+ \sigma f=I.$
Moreover,   $\sigma\sigma=0$  and thus   
\begin{equation}\label{tva}
\sigma(\dbar\sigma)=(\dbar\sigma)\sigma.
\end{equation}
Since $\sigma$ is odd, %%%In view of \eqref{gabbe},
$
\nabla_\End\sigma =\nabla\circ\sigma+\sigma\circ\nabla=
f\sigma+\sigma f-(\dbar\circ\sigma+\sigma\circ\dbar),
$
so we get
\begin{equation}\label{tre}
\nabla_\End \sigma=I-\dbar\sigma.
\end{equation}
%%%\smallskip 
Notice that $\dbar\sigma$ has even degree.
In $X\setminus Z$ we define the $\End E$-valued form, cf., \eqref{tre}, 
\begin{equation}\label{udef}
u=\sigma (\nabla_{\End}\sigma)^{-1}=
\sigma(I-\dbar\sigma)^{-1}=\sigma+\sigma(\dbar\sigma)+\sigma(\dbar\sigma)^2+\ldots.
\end{equation}
Now, 
$
\nabla_\End u=\nabla_\End\sigma(\nabla_\End\sigma)^{-1}
-\sigma \nabla_\End(\nabla_\End\sigma)^{-1},
$
and since $\nabla_\End^2=0$ we thus have
\begin{equation}\label{idol}
\nabla_{\End} u=I.
\end{equation}
Notice that
$$
u=\sum_{\ell}\sum_{k\ge \ell+1}u^\ell_k%%
$$
where 
$$
u^{\ell}_{k}=\sigma_{k}(\dbar\sigma_{k-1})\cdots(\dbar\sigma_{\ell+1})
$$
is in $\E_{0,k-\ell-1}(\Hom(E_\ell,E_{k}))$ over $X\setminus Z$.
In view of \eqref{tva} we also have
\begin{equation}\label{variant}
u^{\ell}_{k}= (\dbar\sigma_{k})(\dbar\sigma_{k-1})\cdots
(\dbar\sigma_{\ell+2})\sigma_{\ell+1}.
\end{equation}
Let
$$
u^\ell=\sum_{k\ge \ell+1}u^\ell_k,
$$
be $u$ composed with the projection $E\to E_\ell$.
We can make a current extension of $u$ across  $Z$ following
\cite{PTY} and the proof of Theorem~1.1 in \cite{A2}. 
In fact,  
after a sequence of suitable resolutions we may assume that
the  sections  $F_j=\det^{\rho_j}f_j$ of 
 $\Lambda^{\rho_j}E_j^*\otimes\Lambda^{\rho_j} E_{j-1}$ are of the form
$F_j=F_j^0F_j'$, where $F_j^0$ is a monomial  and $F_j'$ are  non-vanishing.
If $F$ is a holomorphic function that vanishes on $Z$,  in the same way we may
assume that $F=F^0F'$.
By Lemma~\ref{sing},  $\sigma_j=\alpha_j/F_j^0$, where $\alpha_j$ is
smooth across $Z$.  
Since  $\alpha_{j+1}\alpha_j=0$ outside the set $\{F_{j+1}^0F_j^0=0\}$,
thus    $\alpha_{j+1}\alpha_j=0$   everywhere.  Therefore, cf., \eqref{variant},
it is easy to see that
\begin{equation}\label{polly}
u^\ell_{\ell+k}=\frac{(\dbar\alpha_{\ell+k})(\dbar\alpha_{\ell+k-1})\cdots(\dbar\alpha_{\ell+2})
\alpha_{\ell+1}}
{F^0_{\ell+k}\cdots F^0_{\ell+1}}.
\end{equation}
Since  $F_j$ only vanish on $Z$ and $F$ vanishes there, 
$F^0$ must contain each coordinate factor that occurs in any  $F_j^0$.
It follows now that $\lambda\mapsto|F|^{2\lambda} u$ has a current-valued analytic
continuation to $\Re\lambda>-\epsilon$, and 
that $U=|F|^{2\lambda}u|_{\lambda=0}$ is a current extension of $u$.

In the same way we can now define the residue current $R=R(E_\bullet)$ associated to
\eqref{complex} as
$$
R=\dbar|F|^{2\lambda}\w u|_{\lambda=0}.
$$
It clearly has its support on  $Z$. If 
$
R^\ell_k=\dbar|F|^{2\lambda}\w u^\ell_k|_{\lambda=0}
$
and $R^\ell$ is defined analogously, then
$$
R=\sum_\ell R^\ell=\sum_\ell\sum_{k\ge\ell+1} R^\ell_k.
$$
Notice that  $R^\ell_k$ is a $\Hom(E_\ell,E_k)$-valued $(0,k-\ell)$-current. 
The currents $U^\ell$ and $U^\ell_k$ are defined analogously. 
Notice that $U$ has odd degree and $R$ has even degree.
In analogy with   Theorems~1.1 and 1.2 in \cite{A2}
we  have:

\begin{prop}\label{huvud}
If  $U$ and $R$ are  the currents associated to the complex
\eqref{complex} then
\begin{equation}\label{urformel}
\nabla_{\End}U= I-R, \qquad  \nabla_{\End} R=0.
\end{equation}
Moreover, 
$R^\ell_{k}$ vanishes if $k-\ell<\codim Z$,
and $\bar \xi R=d\bar\xi\w R=0$ if $\xi$ is holomorphic and vanishes on  $Z$. 
\end{prop}

The residue current $R=R(E_\bullet)$ is related to the (lack of) exactness of
the sheaf complex associated to \eqref{complex} in the following way.

\begin{prop}\label{korr}
Let  $R=R(E_\bullet)$ be the residue current associated with  \eqref{complex}  and let
 $\phi$ be a holomorphic section of $E_\ell$.

\noindent (i)  If  
$f_\ell\phi=0$ and $R^\ell\phi=0$, then
locally there is a holomorphic section $\psi$ of $E_{\ell+1}$ such that
$f_{\ell+1}\psi=\phi$.

\noindent  (ii)   If moreover $R^{\ell+1}=0$, then  the existence of such a local
solution $\psi$ implies that  $R^\ell\phi=0$.
\end{prop}

\begin{proof}
Let $U$ be the associated current such that 
\eqref{urformel} holds. Then 
$\nabla (U\phi)=\phi-U(\nabla\phi)  -R\phi$. Since $U\phi= U^\ell\phi$,
$R\phi=R^\ell\phi$, and $\nabla\phi=f_\ell\phi-\dbar\phi$, 
it follows from the  assumptions of $\phi$ that 
$\nabla (U^\ell\phi)=\phi$. 
Now (i) follows by solving a sequence of $\dbar$-equations locally.
%%%
For the second part, assume that $f_{\ell+1}\psi=\phi$. Then by \eqref{urformel},
$
R^\ell\phi= R\phi=R(\nabla\psi)=\nabla(R\psi)=\nabla (R^{\ell+1}\psi)=0.
$
\end{proof}

If now \eqref{plex} is a generically exact holomorphic complex of Hermitian bundles,
since  $\rank f_1$ is generically constant, we can define 
$\sigma_1$ in an unambiguous way in $X\setminus Z$, and therefore
the currents $R^\ell$ for $\ell\ge 0$ can be defined as above, and we have:

\begin{cor}\label{korrkort}
If $R=R(E_\bullet)$ is the residue current associated to
\eqref{plex}, then   Proposition~\ref{korr}
holds (for $\ell\ge 0$), provided that $f_0\phi=0$ is interpreted as
$\phi$ belonging  generically (outside $Z$) to  the image of $f_1$.
\end{cor}

If $f_1$ is generically surjective, in particular  if
$\rank E_0=1$ and $f_1$ is not identically $0$, then this latter condition is
of course automatically fulfilled.

\begin{proof}
The corollary actually follows just from 
a careful  inspection of the arguments in the proof of Proposition~\ref{korr}.
Another way is to  extend \eqref{plex}
to a generically exact complex \eqref{complex} and then refer directly
to Proposition~\ref{korr}, noting that the definition of $R^\ell$ for
$\ell\ge 0$ as well as the condition $f_0\phi=0$ are  independent of 
such an  extension.
\end{proof}

\section{Residue currents with prescribed annihilators}

The exactness of  \eqref{splex2} is characterized by the 
current $R$ associated with \eqref{plex}.

\begin{thm}\label{main1}
Assume that \eqref{plex} is generically exact, let
$R$ be the associated residue current, and let \eqref{splex2}
be the associated complex of sheaves. Then
$R^\ell=0$ for all $\ell\ge 1$ if and only if
\eqref{splex2} is exact.
\end{thm}

For the proof we will use the following characterization of exactness due to
 Buchsbaum-Eisenbud, see \cite{Eis1} Theorem~20.9:
The complex \eqref{splex2} is exact if and only if
\begin{equation}\label{skrot}
\codim Z_j\ge j
\end{equation}
for all $j$,  where, cf., \eqref{rhodef},
$$
Z_j=\{z;\   \rank f_j< \rho_j\}.
$$

\begin{remark}\label{sliskig}
To be precise we will only use the ``only if''-direction. The other direction 
is  actually a consequence of  %%%%Proposition~\ref{huvud}, 
Corollary~\ref{korrkort} and  (the proof of) Theorem~\ref{main1}.
\end{remark}

\begin{proof}[Proof]
From Corollary~\ref{korrkort}  it follows that
\eqref{splex2} is exact if $R^\ell=0$ for $\ell\ge 1$.
%%%
For the converse, let us now assume that
\eqref{splex2} is exact; by the Buchsbaum-Eisenbud theorem then  \eqref{skrot} holds.
We will prove that $R^1=0$; the case when  $\ell>1$ is handled in the same way.
The idea in the proof 
is based on the somewhat vague principle  that a residue current
of bidegree $(0,q)$ cannot be supported on a variety of codimension $q+1$.
Taking this for granted, we notice to begin with that
$R^1_2=\dbar|F|^{2\lambda}\w  \sigma_2|_{\lambda=0}$ is a $(0,1)$-current
and has its support on $Z_2$, which has  codimension at least $2$. 
Hence  $R^1_2$ must vanish according to the vague principle.
Now, $\sigma_3$ is smooth outside $Z_3$, and hence
$R^1_3=\dbar\sigma_3\w R^1_2=0$ outside $Z_3$; thus $R^1_3$ is supported
on $Z_3$ and again, by the same principle,  $R^1_3$ must vanish etc.
To make this into a strict argument we will use  the following simple lemma
which follows from a Taylor expansion.

\begin{lma}\label{yellow}
Suppose that $\gamma(s,\tau)$ is smooth in  $\C\times\C^r$
and that  moreover
$\gamma(s,\tau)/\bar s$ is smooth where  $\tau_1\cdots\tau_k\neq0$.
Then $\gamma(s,\tau)/\bar s$ is smooth everywhere.
\end{lma}

%%\begin{proof}
%%The assumption means that $\gamma(s,\tau)=\bar s\omega(s,\tau)$ where
%%$\tau_1\cdots\tau_r\neq0$ and $\omega$ is smooth outside 
%%$\tau_1\cdots\tau_r=0$.
%%It follows that, for each $\ell$,  $(\partial^\ell/\partial s^\ell)\gamma(0,\tau)=0$ where 
%%$\tau_1\cdots\tau_r\neq0$,
%%and hence by continuity it holds also when $\tau_1\cdots\tau_r=0$.
%%It now follows from a Taylor expansion in $s$ that
%%$\gamma(s,\tau)/\bar s$ is smooth.
%%\end{proof}

After a sequence of resolutions of singularities the action of $R^1_k$ on 
a test form $\xi$ is a finite sum of integrals of the form
$$
\int\dbar|F^0|^{2\lambda}\w
\frac{(\dbar\alpha_k)(\dbar\alpha_{k-1})\cdots(\dbar \alpha_3)\alpha_2}
{F_k^0 F_{k-1}^0\cdots F_3^0 F_2^0}\w
~\tilde\xi\Big|_{\lambda=0}
$$
where $F^0$, $F^0_i$ and $\alpha_i$ are as \eqref{polly} above, 
and where $\tilde\xi$ is the pullback of $\xi$. To be 
precise, there are also cutoff functions involved that we suppress for simplicity.
Observe that $\dbar|F^0|^{2\lambda}$ is a finite sum of terms like
$a\lambda|F^0|^{2\lambda}d \bar s/\bar s$, where $a$ is a positive integer 
and $s$ is just one of the coordinate functions that divide $F^0$.
We need to show that all the corresponding integrals vanish when $\lambda=0$, and to this end
it is enough to show, see, e.g., Lemma~2.1 in \cite{A2},  that
$$
\eta= \frac{d\bar s}{\bar s}\w(\dbar\alpha_k)(\dbar\alpha_{k-1})\cdots(\dbar \alpha_3)\alpha_2\w\tilde\xi
$$
is smooth ($(d\bar s/\bar s)\w \beta$ being smooth for a smooth  $\beta$,
means  that each term of $\beta$ contains a factor $\bar s$ or $d\bar s$). 

Let $\ell$ be the largest index among $2,\ldots,k$ such that $s$ is a factor in $F_\ell^0$ 
(possibly there is no such index at all; then
$\ell$ below is to be interpreted as $1$) 
and let $\tau_1,\ldots,\tau_r$ denote the coordinates 
that divide $F_k^0\cdots F_{\ell+1}^0$. We claim that,  outside 
$\tau_1\cdots\tau_r=0$, the form
\[
\frac{d\bar s}{\bar s}\w\frac{(\dbar\alpha_k)\cdots(\dbar 
\alpha_{\ell+1})}
{F_k^0 \cdots F_{\ell+1}^0}
\w ~\tilde\xi
\]
is smooth.  
This follows by standard arguments, see, e.g., the proof of Lemma 2.2 in 
\cite{PTY} or the proof of Theorem 1.1 in \cite{A2}; in fact, outside 
$Z_k\cap\ldots\cap Z_{\ell+1}$ the $(n,n-\ell+1)$-form 
$(\dbar\sigma_k)\ldots(\dbar\sigma_{\ell+1})\w\xi$ is smooth and 
it must vanish on $Z_\ell$ for degree reasons, since
$Z_{\ell}$ has codimension at least $\ell$.  
Thus the form
\[
\tilde\eta=\frac{d\bar s}{\bar s}\w(\dbar\alpha_k)\cdots(\dbar \alpha_{\ell+1})
\w ~\tilde\xi
\]
is smooth outside $\tau_1\cdots\tau_r=0$. By Lemma \ref{yellow},
applied  to  
$$
\gamma=d\bar s\w (\dbar\alpha_k)\cdots(\dbar \alpha_{\ell+1})\w\tilde\xi,
$$
$\tilde\eta$  is smooth everywhere, and therefore 
$\eta$ is smooth. %%%%the corresponding $s$-integral vanishes 
\end{proof}

If \eqref{splex2} is exact, then,
with no ambiguity, we can  write  $R_k$ rather than  $R^0_k$.

\begin{proof}[Proof of Theorem~\ref{skalle}]
Since a  free resolution of a free sheaf is pointwise exact,
it follows that $Z_N\subset\cdots \subset Z_1=Z$.
Therefore  $u^0$ is smooth outside $Z$ and thus the support of $R$ must
be contained in $Z$. By Theorem~\ref{main1},  $R^1=0$, and so  
the second assertion, the Noetherian property of
$R=R^0$,  follows from  Corollary~\ref{korrkort}.
\end{proof}

Given any coherent sheaf  $\F$  in a Stein manifold $X$ and  compact
subset $K\subset X$, one can always find a resolution 
\begin{equation}\label{bust}
\cdots \to\Ok^{\oplus r_2}\to\Ok^{\oplus r_1}\to \Ok^{\oplus r_0}
\end{equation}
of $\F$ in a neighborhood of $K$, e.g., by
iterated use of Theorem~7.2.1 in \cite{Hor}.
The key stone  in the proof of Theorem~\ref{main1}, the Buchsbaum-Eisenbud  theorem, 
in general requires that the resolution \eqref{bust} starts with $0$ somewhere on the left.
However,
by the Syzygy theorem and Oka's lemma, $\Ker (\Ok^{\oplus r_{\ell}}\to\Ok^{\oplus r_{\ell-1}})$ is
(locally) free for large $\ell$, so we can replace such a module   $\Ok^{\oplus r_\ell}$ with
this kernel and $0$ before that. Therefore Theorem~\ref{main1} holds and we
have

\begin{prop}\label{gatu}
Let  $\J$ be a coherent subsheaf of $\Ok^{\oplus r_0}$ in  a Stein manifold $X$.
For each compact subset $K\subset X$ there is a 
 residue  current $R$ defined in a neighborhood
of $K$ such that $\ann R=\J$.
\end{prop}

%%Notice that in this case $R=(R_k)$, where 
%%$R_k$ is an $r_{k}\times r_0$-matrix  of scalar-valued residue
%%currents. 

%%If $\phi$ is an  $r_0$-column  of functions
%%in $\Ok(K)$  then $R_k\phi$ is an $r_k$-column of currents in a neighborhood
%%of $K$. We can also choose a matrix $f_0$ such that $\phi$ is generically in
%%the image of $f_1$ if and only if $f_0\phi=0$ and we have, 
%%cf., the proofs  of Theorem~\ref{korr} and Corollary~\ref{korrkort}: 

%%\smallskip
%%{\it  The column  $\phi\in\Ok(K)^{r_0}$  of holomorphic functions is in the image of 
%%$\Ok(K)^{r_1}\to\Ok(K)^{r_0}$  if and only if
%%$f_0\phi=0$  and all the residue currents
%%$R_k\phi$ vanish.}

%%%
The  degree of explicitness  of the Noetherian residue  current $R$ in Theorem~\ref{skalle}
is of course directly depending on the degree of explicitness of the  resolution.

\begin{ex}[The Koszul complex]\label{ex1}
Let $H$ be a Hermitian bundle over $X$ of rank $m$ and let $h$ be a non-trivial 
holomorphic section of
the dual bundle $H^*$. Then $h$ can be considered as a morphism $H\to\C\times X$, and
we get a generically exact complex \eqref{plex} by taking
$E_k=\Lambda^k H$ and let all the mappings $f_k$ be interior multiplication
with $f$. If $\eta$ is the section of $E$ over $X\setminus Z$ of minimal norm such that
$f\cdot \eta=1$, then  $\sigma_k\xi=\eta\w\xi$ for
sections $\xi$ of $E_{k-1}$, and hence 
$u^\ell_k=\eta\w(\dbar\eta)^{k-\ell-1},$
acting on $\Lambda^\ell H$ via  wedge multiplication.
Thus 
$R^\ell_k=\dbar|h|^{2\lambda}\w\xi\w(\dbar\xi)^{k-\ell-1}|_{\lambda=0}$
are  precisely the currents  considered in \cite{A2}.
If $h$ is a complete intersection and  $h=h_1 e_1^*+\cdots
+h_me_m^*$  in some local holomorphic frame $e_j^*$ for $H^*$,
then $R$ is precisely  the Coleff-Herrera product  \eqref{ch} 
times $e_1\w\ldots \w e_m$, where $e_j$ is the dual frame, see \cite{A2}.
\end{ex}

%%In \cite{??} and \cite{??} the corresponding residues 
%%obtained from the 
%%The Eagon-Northcott and Buchsbaum-Rim complexes
%%are used. For blabaa see ArXiv version

%%\begin{ex}[The Eagon-Northcott and Buchsbaum-Rim complexes]\label{ex2}
%%Suppose that $H$ and $Q$ are Hermitian bundles of
%%ranks $m$ and $r$ respectively, and $h\colon H\to Q$  is a  generically surjective holomorphic morphism.
%%We then have a natural morphism 
%%$$
%%\det h\colon\Lambda^r H\to \det Q,
%%$$
%%cf., Section~\ref{prelim}.
%%
%%\end{ex}

We now consider  a simple %%
example of a non-complete intersection ideal.

\begin{ex}\label{enklaste}
Consider the ideal $J=(z_1^2,z_1 z_2)$ in $\C^2$  with zero variety $\{z_1=0\}$.
It is easy to see that
\begin{equation}\label{voj}
0\to \Ok\stackrel{f_2}{\longrightarrow}\Ok^{\oplus 2} \stackrel{f_1}{\longrightarrow} \Ok,
\end{equation}
where
$$
f_1=
\left [ \begin{array}{cc}
z_1^2 & z_1z_2
\end{array}\right ] \ \text{ and }\
f_2=
\left [ \begin{array}{c}
z_2 \\ -z_1
\end{array}\right ],
$$
is a (minimal) resolution of $\Ok/J$.
We equip  the corresponding vector bundles with the trivial Hermitian metrics.
Since  $Z$ has codimension  1,  $R$  consists of the two parts
$R_2=\dbar |F|^{2\lambda}\w u_2^0|_{\lambda=0}$ and
$R_1=\dbar |F|^{2\lambda} \w u_1^0|_{\lambda=0}$, 
where
$u_2^0=\sigma_2 \dbar \sigma_1$ and $u_1^0=\sigma_1$, respectively.
%%Notice that $\sigma_1=f_1^*(f_1f_1^*)^{-1}$ and $\sigma_2=(f_2^*f_2)^{-1}f_2^*$.
%%%
To compute
$R$   it is enough to make a simple blow-up at the origin, 
and one gets, cf., \cite{W} and \cite{W2},  that 
$$
R_2=\dbar \bigg[\frac{1}{z_1^2}\bigg ]\wedge
\dbar \bigg[\frac{1}{z_2}\bigg ] \quad \text{and}\quad
R_1=
\left [ \begin{array}{c}
0 \\ 1
\end{array}\right]
\Big[\frac{1}{z_2}\Big]\dbar\Big[\frac{1}{z_1}\Big].
$$
We see that $\ann R_2=(z_1^2,z_2)$ and $\ann R_1=(z_1)$, and hence
$\ann R=(z_1^2, z_2)\cap (z_1) =J$ as expected.
Notice that the Koszul complex associated with the ideal $J$ is
like \eqref{voj} but with an extra factor $z_1$ in 
the mapping $f_2$.
Then the  current $R^0_1$ is of course the same as before, but
$$
R_2^0=\frac{1}{2} \dbar\Big[\frac{1}{z_1^3}\Big]\w\dbar\Big[\frac{1}{z_2}\Big].
$$
In this case $\ann R^0=\ann R^0_2\cap\ann R^0_1=
(z_1^3, z_2)\cap (z_1)$ which is strictly smaller than $J$.
%%
%%%
Roughly speaking,  the annihilator of  $R^0_2$ is too small, since
the singularity of $\sigma_2$ and hence of $u_2^0$  is too big, due to the extra  factor
$z_1$ in $f_2$.
\end{ex}

There has recently been a lot of work done on finding free resolutions of 
monomial ideals, see for example \cite{MS}, \cite{BPS} or \cite{BaS}. 
For more involved explicit computations of residue currents
for monomial ideals, see \cite{W}.
We conclude with a simple  example where $\ann(\Ok(E_0)/J)=0$.

\begin{ex}
Consider the submodule $J$ of $\Ok^{\oplus 2}$ 
generated by  
$f_1= [z_1z_2 \ \   -z_1^2]^T$
and the resolution $0\to \Ok\stackrel{f_1}{\longrightarrow}\Ok^{\oplus 2}$,
which is easily seen to be minimal.
Notice that  $Z=\{z_1=0\}$ is the associated set 
where $\Ok^{\oplus 2}/J$ is not locally free,
or equivalently where $f_1$ is not locally constant.
Moreover,  notice that $\ann(\Ok^{\oplus 2}/J)=0$.
The associated  residue current is  
\[
R=R_1=
\Big[\frac{1}{z_2}\Big]\dbar\Big[\frac{1}{z_1}\Big]
\left [ \begin{array}{cc}
0 & 1
\end{array}\right].
\]
If we extend  the complex  with the mapping 
$f_0=[z_1\  z_1]$
the new complex is still exact outside $Z$. 
%%%
Observe that $\ann R$ is generated by
$ z_1 [1\  1 ]^T$
and moreover that $\Ker f_0$ is generated by
$[ z_2\  -z_1]^T$.
Thus $\Ker f_0\cap \ann R=J$  as expected.
\end{ex}

\section{Cohen-Macaulay ideals and modules}\label{cmsection}

Let $\F_x$ be a  $\Ok_x^r$-module. The minimal length $\nu_x$ of a 
resolution of $\F_x$ is precisely $n-\depth \F_x$,
and  $\depth\F_x\le \dim \F_x$, so
the length of the resolution is at least equal to 
$\codim \F_x$.  Recall that the 
$\F_x$ is Cohen-Macaulay if $\depth\F_x=\dim \F_x$,
or equivalently, $\nu_x=\codim\F_x$,
see \cite{Eis1}.
As usual we say that an ideal $J_x\subset\Ok_x$ is  Cohen-Macaulay
if $\F_x=\Ok_x/J_x$ is a Cohen-Macaulay module.

A  coherent analytic  sheaf $\F$ is 
Cohen-Macaulay if $\F_x$ is Cohen-Macaulay for each $x$.
If we have any locally free resolution of $\F$ and $\codim\F=p$, then 
at each point $\Ker (\Ok(E_{p-1})\to\Ok(E_{p-2}))$
is free by the uniqueness theorem, see below,  so by Oka's lemma the kernel
 is locally free;  hence we can modify the given resolution
to a locally free resolution of minimal length $p$.
Notice that  the   residue current 
associated with a resolution of minimal length $p$ just consists
of the single term $R=R_p^0$, which locally is a 
$r_p\times r_0$-matrix of currents. 

\begin{thm}\label{cmthm1}
Suppose that  $\F$ is a  coherent  analytic sheaf
with codimension $p>0$ that is Cohen-Macaulay, and assume that
\begin{equation}\label{fixp}
0\to\Ok(E_p)\to\cdots\to \Ok(E_1)\to \Ok(E_0)
\end{equation}
is a locally free resolution of $\F$ of minimal length $p$.
Then the associated Noetherian current 
is independent of the Hermitian metric.
\end{thm}

%%Notice that since $p>0$, i.e, $\ann \F \neq 0$,
%%the right-most mapping in
%%\eqref{fixp} is pointwise  surjective outside $Z$.

\begin{proof}
Assume that $u$ and $u'$ are the forms in $X\setminus Z$ constructed by means of
two different choices of metrics on $E$. Then
$\nabla_{\End}u=I$ and $\nabla_{\End}u'=I$ in $X\setminus Z$, and hence
%%if $w=u u'$ we have 
$$
\nabla_{\End}(u u')=(\nabla_{\End}u) u'-u\nabla_{\End}u'=u'-u,
$$
where the minus sign occurs since $u$ has odd order.
For large $\Re\lambda$  we thus have,
cf., the proof of Proposition~\ref{huvud},
$$
\nabla_{\End}\big(|F|^{2\lambda} uu'\big)=|F|^{2\lambda}u'-|F|^{2\lambda}u-\dbar|F|^{2\lambda}\w uu'.
$$
As before one can verify  that  each term
admits an  analytic continuation to $\Re\lambda>-\epsilon$, and 
 evaluating at $\lambda=0$ we get
$
\nabla_{\End} W=U'-U-M,
$
where $W=|F|^{2\lambda}uu'|_{\lambda=0}$,
and   $M$ is the residue current 
\begin{equation}\label{mdef}
M=\dbar|F|^{2\lambda}\w uu'|_{\lambda=0}.
\end{equation}
Since $\nabla_{\End}^2=0$, by Proposition~\ref{huvud} we therefore get
\begin{equation}\label{gurka}
R-R'=\nabla_{\End} M.
\end{equation}
However, since the complex ends up at $p$, each term in 
$uu'$ has at most bidegree $(0,p-2)$ and hence
the current $M$ has at most bidegree $(0,p-1)$. Since it is
supported on $Z$ with codimension $p$,  it must vanish, cf., the proof
of Proposition~\ref{huvud}.
\end{proof}

When  $\F=\Ok(E_0)/\J$ is Cohen-Macaulay 
we can also define a cohomological residue  that characterizes
the module sheaf $\J=\Im(\Ok(E_1)\to\Ok(E_0))$ locally.
Suppose that  we have a fixed resolution \eqref{fixp} of minimal length 
and let us assume that $p>1$.
If  $u$  is  any solution to $\nabla_{\End} u=I$ in $X\setminus Z$,
then   $u^0_p$ is a $\dbar$-closed
$\Hom(E_0, E_p)$-valued $(0,p-1)$-form. Moreover if $u'$ is another
solution, then it follows from the preceding proof that
$\dbar(uu')^0_p=u^0_p-u'^0_p$.  Therefore $u_p^0$ defines a 
Dolbeault cohomology class $\omega\in H^{0,p-1}(X\setminus Z,\Hom(E_0,E_p))$.
If $\phi$ is a holomorphic section of $E_0$ then
$\omega\phi=[u^0_p\phi]$ is an element in $H^{0,p-1}(X\setminus Z,E_p)$.
Moreover,   if $v$ is any solution in $X\setminus Z$ to
$\nabla v=\phi$, then $v_p$ defines the class $\omega\phi$.
In fact, $\nabla (u v)=v-u\phi=v-u^0\phi$ so that
$\dbar (uv)_p =u^0_p\phi -v_p$.
Precisely as for a complete intersection, \cite{DS} and \cite{P1},
we have the following  cohomological duality principle.

\begin{thm}\label{dualthm}
Let $X$ be a Stein manifold and let \eqref{fixp} be a resolution of minimal
length $p$ of the Cohen-Macaulay sheaf $\Ok(E_0)/\J$ over $X$, and assume that $p>1$. 
Moreover, let
$\omega$ be the associated  class in $H^{0,p-1}(X\setminus Z,\Hom(E_0,E_p))$. 
For a  holomorphic section $\phi$ of $E_0$ the following conditions are equivalent:

\smallskip

(i) $\phi$ is a global section of $\J$.

\smallskip
(ii) The class $\omega\phi$ in $X\setminus Z$ vanishes.

\smallskip

(iii)  $\int \omega\phi\w\dbar\xi=0$ for all $\xi\in\D_{n,n-p}(X, E_p^*)$
such that $\dbar\xi=0$ in a neighborhood of $Z$.
\end{thm}

Notice that if $R$ is the associated  Noetherian current, then $\dbar U^0_p=R_p$,
so by Stokes'  theorem, $(iii)$ is equivalent to  that
$\int R_p\phi\w\xi=0$ for  all $\xi\in\D_{n,n-p}(X, E_p^*)$
such that $\dbar\xi=0$ in a neighborhood of $Z$.

\smallskip
If $p=1$, then $f_1$ is an isomorphism outside $Z$, so its inverse
$\omega=\sigma_1$ is a holomorphic $(0,0)$-form in $X\setminus Z$.
Thus a holomorphic section $\phi$ of $E_0$ belongs to  $\J$ if and only if
$\omega\phi$ has a holomorphic extension across $Z$.

\begin{proof}%%[Sketch of proof]
If $(i)$ holds, then $\phi=f_1\psi$ for some holomorphic $\psi$;
thus $\nabla \psi=\phi$. However, since $p>1$, $\psi$ has no
component in $E_p$, and hence by definition  the class  $\omega\phi$ vanishes.
The implication $(ii)\to(iii)$ follows from Stokes' theorem.  

Let us now assume that $(iii)$ holds, and choose a point $x$ on $Z$.
Let $v_k= u_k^0\phi$.
If $X'$ is an appropriate small neighborhood of $x$, then,
since $Z$ has codimension $p$  and $v_p$ is a $\dbar$-closed $(0,p)$-current,
one can verify that the condition $(iii)$  ensures that 
$\dbar w_p=v_p$ has a solution 
in  $X'\setminus \overline W$, where $W$ is a  small
neighborhood of $Z$ in $X'$. Then, successively,
 all the lower degree  equations
$\dbar w_k=v_k+f_{k+1}w_{k+1}$, $k\ge 2$,  can be solved 
in similar domains. Finally, we get
a holomorphic solution $\psi=v_1+f_2 w_2$ to $f_1\psi=\phi$,
in such a domain. By Hartogs' theorem $\psi$ extends across $Z$ in $X'$.
Alternatively, one can obtain such a local holomorphic solution $\psi$, 
using the decomposition  formula \eqref{decomp} below
and  mimicking the proof of the corresponding statement for a complete
intersection in \cite{P1}; cf., also the proof of  Proposition~7.1 in \cite{A7}.
Since $X$ is Stein, one can piece together
to a global holomorphic solution to $f_1\psi=\phi$, and hence   $\phi$
is a section of $\J$.
\end{proof}

\begin{ex} 
Let $J$ be an ideal in $\Ok_0$ of dimension zero. Then it is Cohen-Macaulay and 
for each germ  $\phi$  in $\Ok_0$,  $\omega\phi$ defines a
functional on $\Ok_0(E_n^*)\simeq\Ok_0^{r_n}$.
If $J$ is defined by a complete intersection, then 
we may assume that \eqref{fixp} is the Koszul complex. Then
$r_n=1$, and in view of the Dolbeault isomorphism,
see, e.g., Proposition~3.2.1 in \cite{P1},
$\omega\phi$ is just the classical Grothendieck residue.
\end{ex}

\smallskip

For the rest of this section we will restrict our attention
to modules over the local ring  $\Ok_0$, and we let
$\Ok(E_k)$ denote the free $\Ok_0$-module of germs of
holomorphic sections at $0$ of the vector bundle $E_k$.
%%ablbala analogous for ring $\Ok(K)$  $K$ semi-analytic balbala.
%%
Given a  free resolution \eqref{splex2}
%%\begin{equation}\label{fres}
%%0\to \Ok(E_N)\stackrel{f_p}{\longrightarrow}\cdots
%%\stackrel{f_2}{\longrightarrow}\Ok(E_1)\stackrel{f_1}{\longrightarrow}\Ok(E_0)
%%\end{equation}
of a module $\F_0$  over $\Ok_0$ and given metrics on $E_k$ we thus get
a germ $R$ of a Noetherian residue current at $0$. 
Recall that the resolution \eqref{splex2}
is {\it minimal}  if for each $k$,  $f_k$ maps  a basis of 
$\Ok(E_k)$ to a minimal set of generators of $\Im f_k$.
The uniqueness theorem, see, e.g., Theorem~20.2 in \cite{Eis1},
states that  any two minimal (free) resolutions are equivalent, and moreover,
that any (free) resolution has a minimal resolution as a direct summand.

For a Cohen-Macaulay module $\F_0$ over $\Ok_0$  we  
have the following uniqueness.

\begin{prop}
Let $\F_0$ be a Cohen-Macaulay module over $\Ok_0$ 
of  codimension $p$. If we have two minimal
free resolutions $\Ok(E_\bullet)$ and $\Ok(E_\bullet')$ of
$\F_0$, then there are holomorphic invertible matrices
$g_p$ and $g_0$ (local holomorphic isomorphism $g_p\colon E_p'\simeq E_p$
and $g_0\colon E'_0\simeq E_0$)  such that $R=g_p R'g_0^{-1}$.
\end{prop}

Since minimal resolutions have minimal length $p$, the currents are independent
of the metrics, in view of Proposition~\ref{cmthm1}.

\begin{proof} By the uniqueness theorem there are 
holomorphic local isomorphisms $g_k\colon E_k'\to E_k$ such that 
$$
 \begin{array}{ccccccccccc}
0 & \to & \Ok(E_p') & \stackrel{f_p'}{\longrightarrow} &  \cdots & \stackrel{f_2'}{\longrightarrow} & 
\Ok(E_1') & 
\stackrel{f_1'}{\longrightarrow} & \Ok(E_0')\\ 
{} & &    g_p\downarrow &&  &&       g_1\downarrow  &&  {g_0}\downarrow & \\
0 & \to & \Ok(E_p) & \stackrel{f_p}{\longrightarrow} &  
\cdots & \stackrel{f_2}{\longrightarrow} & \Ok(E_1) & 
\stackrel{f_1}{\longrightarrow} & \Ok(E_0)
\end{array}
$$
commutes.   Let $g$ denote the induced
isomorphism $E\to E'$. Choose any metric on $E$ and equip $E'$ with the
induced metric, i.e., such that $|\xi|=|g^{-1}\xi|$ for a section
$\xi$ of $E'$. 
If $\sigma\colon  E\to E$ and $\sigma'\colon E'\to E'$ are the associated
endomorphisms over $X\setminus Z$, cf.,  Section~\ref{hungrig},
then $\sigma' = g\sigma g^{-1}$ in $X\setminus Z$,  and therefore
$$
u'=\sigma' +(\dbar\sigma')\sigma' +\cdots=g(\sigma+(\dbar\sigma)\sigma+\cdots)g^{-1}
=g u g^{-1}.
$$
Therefore,  $(u')^0_p=g_p u^0_pg_0^{-1}$,  and hence  the  statement follows
since $R=R_p=R_p^0$. 
\end{proof}

We shall now consider the residue  current associated to a general
free resolution.

\begin{thm}\label{cmthm2}
Let $\F_0$ be a Cohen-Macaulay module over $\Ok_0$ 
of  codimension $p$.
If $R$ is the  residue current  associated to an
arbitrary free resolution \eqref{splex2} (and given metrics  on $E_k$)
and $R'=R_p'$ is associated to a minimal resolution
$
0\to\Ok(E_p')\stackrel{f_p'}{\longrightarrow}\cdots
\stackrel{f_2'}{\longrightarrow} \Ok(E_1')\stackrel{f_1'}{\longrightarrow}\Ok(E_0'),
$
then 
\begin{equation}\label{gamas}
R_p=h_p R_p'\beta_0,
\end{equation}
where $\beta_0\colon E_0\to E_0'$ is a local holomorphic  pointwise surjective  morphism and
$h_p$ is a  local smooth  pointwise injective  morphism $h_p\colon E_p'\to E_p$.
Moreover, for each $\ell>0$,
$$
R_{p+\ell}=\alpha_\ell R_p,
$$
where $\alpha_\ell$ is a  smooth $\Hom(E_p,E_{p+\ell})$-valued $(0,\ell)$-form.
\end{thm}

%%It is clear from this result that the annihilaa

\begin{proof}
By the uniqueness theorem for resolutions, the resolution $E_\bullet'$ is
isomorphic to a direct summand in $E_\bullet$, and in view of the preceding
proposition, we may assume that
%%It is well-known,  see, e.g., \cite{Eis1}, that  \eqref{fres}
%%has a minimal resolution as a direct summand, 
%%and since all minimal resolutions are equivalent, 
$$
\Ok(E_k)=\Ok(E_k'\oplus E_k'')=\Ok(E_k')\oplus\Ok(E_k'')
$$
and $f_k=f_k'\oplus f_k''$,  so that 
$$
\begin{array}{ccccccccccc}
0 & \to & \Ok(E_p') & \stackrel{f_p'}{\longrightarrow} &  \cdots & \stackrel{f_2'}{\longrightarrow} & 
\Ok(E_1') & 
\stackrel{f_1'}{\longrightarrow} & \Ok(E_0')\\ 
{i_{p+1}\downarrow} & &    i_p\downarrow &&  &&       i_1\downarrow  &&  {i_0}\downarrow & \\
\to \Ok(E_{p+1}) & \stackrel{f_{p+1}}{\longrightarrow} & \Ok(E_p) & \stackrel{f_p}{\longrightarrow} &  
\cdots & \stackrel{f_2}{\longrightarrow} & \Ok(E_1) & 
\stackrel{f_1}{\longrightarrow} & \Ok(E_0),
\end{array}
$$
where $i_k\colon E_k'\to E_k'\oplus E_k''$ are the natural injections, and 
$$
\to \Ok(E_{p+1}'')\stackrel{f_{p+1}''}{\longrightarrow}
\Ok(E_p'')\stackrel{f_p''}{\longrightarrow}\cdots
\stackrel{f_2''}{\longrightarrow} \Ok(E_1'')
\stackrel{f_1''} {\longrightarrow} \Ok(E_0'')
$$
is a resolution of $0$.
In particular,
$$
\to E_{p+1}\stackrel{f_{p+1}''}{\longrightarrow}
E_p''\stackrel{f_p''}{\longrightarrow}\cdots
\stackrel{f_2''}{\longrightarrow} E_1''\stackrel{f_1''}{\longrightarrow} E_0''\to 0
$$
is a pointwise exact sequence of vector bundles, and therefore
the set $Z_k$ where $\rank f_k$ is not optimal coincides with the
set $Z_k'$ where $\rank f_k'$ is not optimal. In particular, $Z_k=\emptyset$
for $k>p$. If we choose, to begin with,
 Hermitian metrics on $E_k$ that respect this direct
sum, and let $\sigma_k$, $\sigma_k'$,  and $\sigma_k''$ be the corresponding
minimal inverses, then 
$
\sigma_k=\sigma_k'\oplus\sigma_k''
$
and hence
$$
u_k^0=(\dbar\sigma_k'\oplus\dbar\sigma_k'')(\dbar\sigma_{k-1}'\oplus\dbar\sigma_{k-1}'')\cdots
(\dbar\sigma_2'\oplus\dbar\sigma_2'')(\sigma_{1}'\oplus \sigma_1'')=
(u')^0_k\oplus(u'')^0_k
$$
for all $k$. However, $(u'')^0_k$ is smooth, and hence
$$
R_p=R_p'\oplus 0,\quad R_k=0\ \text{for} \ k\neq p.
$$
For this particular choice of metric thus
\eqref{gamas} holds with $h_p$ as the natural injection $i_p\colon E_p'\to E_p$
and $\beta_0$ as the natural projection.

Without any risk of confusion we can therefore from now on let $R'_p$ denote the residue current
with  respect to this particular metric on $E$, and moreover let 
$\sigma'$ denote the minimal inverse of $f$ with respect to this metric etc.
We now choose other metrics on $E_k$ and let 
$R_k$ from now on denote the  residue current associated with this new metric.
Following the notation in
the proof of Proposition~\ref{cmthm1} we again have \eqref{gurka},  and for degree 
reasons still $M_p^0=0$; here  $M^\ell_k$ denotes the component of $M$ that
takes values in $\Hom(E_\ell,E_k)$.   Thus
$$
R_p-R'_p=f_{p+1}M^0_{p+1}.
$$
Moreover, if we expand $uu'$, we get
\begin{multline*}
M^0_{p+1}=
\dbar|F|^{2\lambda}\w\big[
\sigma_{p+1} \sigma_p'(\dbar\sigma_{p-1}')\cdots(\dbar\sigma_1')+\\
\sigma_{p+1} (\dbar\sigma_p) \sigma_{p-1}'(\dbar\sigma_{p-2}')\cdots(\dbar\sigma_1')+
\cdots\big]|_{\lambda=0}.
\end{multline*}
However, $\sigma_{p+1}(\dbar\sigma_p)=(\dbar\sigma_{p+1})\sigma_p$
and $\sigma_{p+1}$ is smooth since $Z_{p+1}$ is empty, so 
$$
M^0_{p+1}=-\sigma_{p+1}R'_p+(\dbar\sigma_{p+1})M_p^0=-\sigma_{p+1}R'_p.
$$
Thus,  
$$
R_p=R_p'-f_{p+1}\sigma_{p+1}R_p'=(I_{E_p}-f_{p+1}\sigma_{p+1})R_p'.
$$
Since $f_{p+1}$ has constant rank,  $H=\Im f_{p+1}$ is a smooth subbundle
of $E_p$. 
Notice that 
$\Pi=I_{E_p}-f_{p+1}\sigma_{p+1}$ is the orthogonal projection of
$E_p$ onto the orthogonal complement of $H$ with respect to the new metric.
In this case therefore $h$ in 
\eqref{gamas} becomes  the natural injection $i_p\colon
E_p'\to E_p$ composed by $\Pi$, and 
since $E_p'\cap H=0$,  $h$  is  pointwise injective.

Since $Z_k$ is empty for $k>p$, $\sigma_k$ is smooth for $k>p$ and hence for $\ell>p$,
$$
R_\ell=\dbar|F|^{2\lambda}\w(\dbar\sigma_\ell)\cdots(\dbar\sigma_{p+1})
u^0_p=(\dbar\sigma_\ell)\cdots(\dbar\sigma_{p+1})\dbar|F|^{2\lambda}\w
u^0_p=\alpha_\ell R_p
$$
where
$\alpha_\ell=(\dbar\sigma_\ell)\cdots(\dbar\sigma_{p+1})$.
\end{proof}

\section{Division and interpolation formulas}\label{intf}

To obtain   formulas for division and interpolation
 that involve our currents $R$ and $U$ we
will use the general scheme developed in \cite{A7}.
Let $z$ be a fixed point in $\C^n$, let 
$\delta_{\zeta-z}$ denote interior multiplication by the vector field
$
2\pi i\sum_1^n (\zeta_j-z_j)(\partial/\partial \zeta_j),
$ 
and let $\nabla_{\zeta-z}=\delta_{\zeta-z}-\dbar$.
Let  $g=g_{0,0}+\cdots +g_{n,n}$  be   a smooth form 
such that $\nabla_{\zeta-z}g=0$ and $g_{0,0}(z)=1$
(here lower indices denote bidegree); 
such a form will be called a weight with respect to the point $z$.
If $g$ has compact support then 
\begin{equation}\label{axel}
\phi(z)=\int g\phi
\end{equation}
for   $\phi$  that are  holomorphic in a neighborhood of the support of $g$,
\cite{A7}.

Let  $D$ be   a ball with center at the origin  in $\C^n$ and
 let $\chi$ be  a cutoff function that is $1$ in a  neighborhood of $\overline{D}$.
Then for each $z\in\overline D$, 
\begin{equation}\label{vikt}
g=\chi-\dbar\chi\w\frac{s}{\nabla_{\zeta-z} s}=
\chi-\dbar\chi\w [s+s\w\dbar s +\cdots +s\w(\dbar s)^{n-1}]
\end{equation}
is a weight,  and it  depends holomorphically on $z$.
Assume that  $\eqref{complex}$ is  a complex of (trivial) bundles over
a neighborhood of  $\overline D$  and let
$\J=\Im f_1$. Let us also fix global frames for the bundles $E_k$.
Then $E_k\simeq\C^{\rank E_k}$ and 
the morphisms $f_k$ are just matrices of holomorphic functions.
One can find (see \cite{A7} for explicit choices)
$(k-\ell,0)$-form-valued holomorphic Hefer morphisms, i.e.,  matrices,  
$H^\ell_k\colon E_k\to E_\ell$   depending holomorphically on $z$ and $\zeta$, such that 
$H^\ell_k=0$ for $ k<\ell$, $ H^\ell_\ell=I_{E_\ell}$, and in general,
\begin{equation}\label{Hdef}
\delta_{\zeta-z} H^\ell_{k}=
H^\ell_{k-1} f_k -f_{\ell+1}(z) H^{\ell+1}_{k};
%%%\quad k\ge\ell;
\end{equation}
here $f$ stands for $f(\zeta)$. Let
$$
HU=\sum_\ell H^{\ell+1} U = \sum_{\ell k} H^{\ell+1}_k U^\ell_k,
\quad  HR=\sum_\ell H^\ell R=\sum_{\ell k} H^\ell_k R^\ell_k.
$$
Then
$
g'= f(z)HU + HUf+ HR
$
maps a section of $E_\ell$ depending on $\zeta$ into a (current-valued) 
section of $E_{\ell }$ depending on both $\zeta$ and $z$. 
Moreover,
$\nabla_{\zeta-z} g'=0 \quad \text{and}\quad g'_{0,0}=I_E$.
If $g$ is  weight with compact support,  
cf.,  Proposition~5.4 in \cite{A7},  we therefore have the representation
\begin{equation}\label{decomp}
\phi(z)=f_{k+1}(z)\int_\zeta H^{k+1}U\phi\w g+
\int_\zeta H^kU f_k \phi\w g+
\int_\zeta H^kR\phi\w g, 
\end{equation}
$z\in \overline D$, for $\phi\in\Ok(\overline D, E_k)$.
Thus  we get an explicit  realization (in terms of $U$) 
of to $f_{k+1}\psi=\phi$,
if  $f_k\phi=0$ and $R\phi=0$, and thus
an explicit proof of Proposition~\ref{korr}~(i).

If we have a complex \eqref{plex} over a neighborhood of $\overline D$, and
either $f_1$ is generically surjective or
we have an extension to a generically exact complex ending at
$E_{-1}$, then \eqref{decomp} still holds for $k=0$. 
If $R$ is Noetherian, then   the last two terms vanish 
if and only if $\phi$ is in $\J$.  We thus obtain an explicit realization
of the membership of $\J$.

\smallskip

In the same way as in \cite{A3}   one can extend these formulas slightly,
to obtain   a characterization of the module  $\E J$ of smooth tuples of functions
generated by $J$, i.e., the set of all $\phi=f_1\psi$ for smooth $\psi$.
For simplicity we assume that $\Ok(E_0)/J$ has positive codimension so that $f_0=0$.
Let $R$ be a Noetherian current for $J$.
First notice that if $\phi= f_1\psi$, then, cf., Proposition~\ref{huvud}, 
$R\phi=R^0\phi=R^0f_1\psi-R^1\dbar\psi=R\nabla\psi=\nabla R^1\psi=0$, so that
$R\phi=0$.
Since each partial derivative $\partial/\partial\bar z_j$ commutes with
$f_1$, we get that
\begin{equation}\label{avil}
R(\partial^\alpha\phi/\partial\bar z^\alpha)=0
\end{equation}
for all multiindices $\alpha$.
The converse can be  proved  by integral formulas precisely as in \cite{A3}, and
thus we have

\begin{thm}\label{smoothy}
Assume that $J\subset\Ok^{\oplus r_0}$ is a coherent subsheaf  such that  $\Ok^{\oplus r_0}/J$ has
positive codimension, and let  $R$ be  a Noetherian residue current for $J$.
Then an $r_0$-tuple $\phi\in \E^{\oplus r_0}$ of smooth functions is in $\E J$
if and only if \eqref{avil} holds for all $\alpha$.
\end{thm}

%%One can also obtain analogous results for  lower regularity as in  \cite{A3} and \cite{A7},
%%as well as a version  where the codimension of  $\Ok^{\oplus r_0}/J$ is zero;
%%one then must add the compatibility condition $f_0\phi=0$.

%%In the case of a complete intersection, Bj\"ork, \cite{JEB3},  has recently given
%%a simple proof of Theorem~\ref{smoothy} based on a deep criterion
%%for membership of ideals of smooth functions in terms of formal power series
%%due to Malgrange, \cite{Mal}.  It extends to a general ideal
%%if our current $R$ is replaced by a tuple of Coleff-Herrera currents
%%$\gamma_j$ such that  $I=\cap\ann \gamma_j$.
%%%

%%\begin{remark}
%%One should notice that  the corresponding statement, where ``smooth''
%%is replaced by ``real-analytic''  easily follows from the holomorphic case.
%%In fact, if $\phi(\zeta)$ is real-analytic, then
%%$\phi(\zeta)=\tilde\phi(\zeta,\bar\zeta)$, where 
%%$$
%%\tilde\phi(\zeta,\omega)=\sum_\alpha\frac{\partial^\alpha\phi}{\partial\bar\zeta^\alpha}(\zeta)
%%(\omega-\bar\zeta)^{\alpha}/\alpha!
%%$$
%%is holomorphic in a neighborhood of $(\zeta,\bar\zeta)$ in $\C^n\times\C^n$.
%%Notice that $R\otimes 1$ is a Noetherian current for $J\otimes 1$ in $\C^n\times\C^n$.
%%If \eqref{avil} holds, it follows that $R\otimes 1\tilde\phi=0$;   hence
%%$f_1(\zeta)\psi(\zeta,\omega)=\tilde\phi(\zeta,\omega)$ and thus 
%%$f_1(\zeta)\psi(\zeta,\bar\zeta)=\phi(\zeta)$.
%%\end{remark}

Let $J$ be a  coherent  Cohen-Macaulay ideal sheaf of codimension $p$
over  some pseudoconvex set $X$ and let $\mu$ be an analytic functional that annihilates $J$.
In \cite{DGSY} was proved (Theorem 4.4) that  $\mu$ can be represented
by an $(n,n)$-current $\tilde\mu$ with compact support of the form
$\tilde\mu=\alpha\w R$, where $\alpha$ is a smooth $(n,n-p)$-form with
compact support and $R$ is the  Coleff-Herrera product
of a complete intersection ideal contained in $J$. In particular,
$\tilde\mu$ vanishes on  $\E J$.
As another  application of our  integral formulas we prove
the following more general result.

\begin{thm}
Let $X$ be a pseudoconvex set in $\C^n$ and let $J$ be a coherent subsheaf of
$\Ok(E_0)\simeq\Ok^{\oplus r_0}$ such that $\Ok(E_0)/J$ has positive codimension.
If $\mu\in\Ok'(X,E_0^*)$ is an analytic functional that vanishes
on $J$, then there is an  $(n,n)$-current $\tilde\mu$ with compact support that
represents $\mu$, i.e.,
\begin{equation}\label{bullar}
\mu.\xi=\tilde\mu.\xi, \quad \xi\in\Ok(X,E_0),
\end{equation}
and such that $\tilde\mu$ vanishes on $\E J$. More precisely we can choose $\tilde\mu$
of the form 
$$
\tilde\mu=\sum_k\alpha_kR_k,
$$
where $R$ is a Noetherian residue current for $J$ and $\alpha_k\in\D_{n,n-k}(X,E_k^*)$.
\end{thm}

Here $E_k$ refers to the trivial vector bundles associated to a free resolution
of $\Ok(E_0)/J$.

\begin{proof}
Assume that $\mu$ is carried by the $\Ok(X)$-convex compact subset $K\subset X$
and let  $V$ be an open neighborhood of $K$.  For each $z\in V$ we can choose a weight
$g^z$ with respect to $z$, such that $z\mapsto g^z$ is holomorphic in $V$ and 
all $g^z$ have support in some compact $\tilde K\subset  X$, see Example~10 in \cite{A2}.
Let $R$ be a  residue current for $J$,
associated to a free resolution of $\Ok(E_0)/J$  in a neighborhood of $\tilde K$, cf,
Proposition~\ref{gatu}.
Now consider the corresponding decomposition \eqref{decomp} (with $k=0$)
that holds for
$z\in V$, with $g=g^z$;  notice that $f_0=0$ by the assumption on  $J$.
The analytic functional $\mu$ has a continuous extension to $\Ok(K,E_0)$ and 
since $\Ok(X)$ is dense in $\Ok(K)$
$\mu$  will vanish on  the first term on the right hand
side in \eqref{decomp}.
If we   define the $(n,n)$-current
$$
\tilde\mu = \mu_z(g^z\w H^0)R=
\sum_k\mu_z(g^z_{n-k,n-k}\w H^0_k) R_k=\sum_k\alpha_k R_k,
$$
then $\alpha_k$ have compact support and  \eqref{bullar} holds.
Since $R$ is Noetherian, $\tilde\mu$ annihilates $\E J$.
\end{proof}

\section{Homogeneous residue currents}\label{rotta}

We will now make a construction of  homogeneous Noetherian residue currents
in $\C^{n+1}$.  This  is  the key to find global Noetherian
currents for polynomial ideals in $\C^n$ by homogenization in the next section.
Let  $S=\C[z_0,z_1,\ldots,z_n]$ be the graded ring of polynomials in
$\C^{n+1}$, and  let  $S(-d)$ be equal to $S$ considered as an $S$-module, but with the
grading shifted by $-d$, so that   the constants
have degree $d$, the linear forms have degree $d+1$ etc.
Assume that 
\begin{equation}\label{gcomplex}
0\to M_N\to\cdots \to M_1\to M_0 
\end{equation}
is  a complex of free graded  $S$-modules,
where 
$$
M_0=S^{\oplus r_0}, \quad 
M_k=S(-d^k_1)\oplus\cdots \oplus S(-d_{r_k}^k).
$$
Then the (degree preserving) mappings are given by matrices of homogeneous elements in $S$.
We can associate to \eqref{gcomplex} a  generically exact
complex of vector bundles 
\eqref{plex} over $\P^n$
%%\begin{equation}\label{ecomplex}
%%0\to E_N\stackrel{f_N}{\longrightarrow}\ldots\stackrel{f_3}{\longrightarrow} 
%%E_2\stackrel{f_2}{\longrightarrow}
%%E_1\stackrel{f_1}{\longrightarrow}E_0,
%%\end{equation}
in the following way. Let $\Ok(\ell)$ be the holomorphic line bundle over $\P^n$ whose sections are 
(naturally identified with)  $\ell$-homogeneous functions in $\C^{n+1}$.
Moreover, let $E^i_j$ be disjoint trivial line bundles over $\P^n$ and let
$$
E_k=\big(E_1^k\otimes \Ok(-d_1^k)\big)\oplus\cdots \oplus \big(E_{r_k}^k\otimes \Ok(-d^k_{r_k})\big).
$$
Notice that  homogeneous elements in $M_\ell$ of degree $r$ precisely corresponds
to the global holomorphic sections of the bundle $E_\ell\otimes\Ok(r)$.

The mappings in \eqref{gcomplex} induce  vector bundle morphisms $f_k\colon E_k\to E_{k-1}$.
We equip $E_k$ with the natural  Hermitian metric, i.e., such that 
$$
|\xi(z)|^2_{E_k}=\sum_{j=1}^{r_k}|\xi_j(z)|^2 |z|^{2d_j^k},
$$
if  $\xi=(\xi_1,\ldots,\xi_{r_k})$, and we  have the associated
currents $U$ and  $R$ as before; they are  associated to the complex
\begin{equation}\label{bula}
0\to E_N\otimes\Ok(r)\stackrel{f_N}{\longrightarrow}\ldots
\stackrel{f_2}{\longrightarrow}
E_1\otimes\Ok(r)\stackrel{f_1}{\longrightarrow}E_0\otimes\Ok(r) 
\end{equation}
as well.

\begin{ex}\label{badlakan}
For each $j,k$ let $\epsilon^k_j$ be a global frame element for the bundle $E^k_j$. 
Then 
$$
R^\ell_k=\sum_{i=1}^{r_\ell}\sum_{j=1}^{r_k} (R^\ell_k)_{ij}\otimes 
\epsilon^k_i\otimes (\epsilon^{\ell}_j)^*,
$$
where  each $(R_k^\ell)_{ij}$ is a $(0,k-\ell)$-current on $\P^n$, taking
values in 
\newline 
$\Hom(\Ok(-d_j^\ell),\Ok(-d_i^k))\simeq\Ok(d_j^\ell-d_i^k)$; 
alternatively $(R_k^\ell)_{ij}$  can be viewed as a 
$(d_j^\ell-d_i^k)$-homogeneous current on $\C^{n+1}\setminus\{0\}$.
In the affine
part $\U_0=\{[z]\in\P^n;\ z_0\neq 0\}$ we have, for each $k$,  a holomorphic frame
$$
e^k_j= z_0^{-d^k_j}\epsilon^k_j, \quad j=1,\ldots, r_k,
$$
for the bundle $E_k$. 
%%Expressed in these local frames, of course  
%%$U$ and $R$ become  residue currents  in $\U_0\simeq\C^n$,
%%satisfying  the relation $\nabla_\End U=I-R$.
In  these   frames
\begin{equation}\label{glunt}
R^\ell_k=\sum_{i=1}^{r_\ell}\sum_{j=1}^{r_k} (\hat R^\ell_k)_{ij}\otimes e^k_i\otimes (e^{\ell}_j)^*,
\end{equation}
where $(\hat R^\ell_k)_{ij}$ are (scalar-valued)  currents  in $\U_0\simeq\C^n$. 
Since $(\hat R^\ell_k)_{ij}$ are    the dehomogenizations of  
$(R^\ell_k)_{ij}$, and $d_j^\ell-d_i^k\le 0$,
it is easily seen  that $(\hat R^\ell_k)_{ij}$ have   current extensions
to $\P^n$. 
\end{ex}

If \eqref{gcomplex} is exact, then according to the Buchsbaum-Eisenbud theorem
for  graded rings, see \cite{Eis2}, the set in $\C^{n+1}$ (or equivalently in $\P^n$)
where the rank of $f_k$ is strictly less than the generic rank $\rho_k$,
has at least codimension $k$. It follows, cf., the proof of Theorem~\ref{main1},  that
$R^\ell=0$ for $\ell\ge 1$, and \eqref{bula} is exact.
In particular, $R=R^0$ is a Noetherian residue
current for the subsheaf $\J\otimes\Ok(r)$ of $\Ok(E_0\otimes\Ok(r))$  generated by $f_1$.
Let now $\phi$ be a global holomorphic section $\phi$ of $E_0\otimes\Ok(r)$, 
that is generically in the image of $f_1$, and such that $R\phi=0$.
Then $\nabla(U^0\phi)=\phi$, cf., the proof of  Proposition~\ref{korr}, and
we obtain a holomorphic section $\psi$ of $E_1\otimes\Ok(r)$ such that
$f_1\psi=\phi$,  provided that we can  solve globally a sequence of $\dbar$-equations.
The first one is $\dbar w_\mu =U^0_\mu\phi$, $\mu=\min(N-1,n)$, and the right hand side
here is a $(0,\mu-1)$-current with values in
$$
E_\mu\otimes\Ok(r)\simeq\oplus_{j=1}^{r_\mu}\Ok(r-d_j^\mu).
$$
Recall that
$H^{0,q}(\P^n, \Ok(\nu))=0$ for all $\nu$ if 
$0<q<n$, whereas $H^{0,n}(\P^n,\Ok(\nu))=0$ if $\nu\ge -n$,
see, e.g.,   \cite{Dem}.
Therefore the equation has a global solution if either $N\le n$ or
$\max_j(r-d_j^n)\ge -n$.
The other equations to solve, $\dbar w_k= U^0_k\phi+f_{k+1}w_{k+1}$,
have lower degree so then there are no
cohomological obstructions. Thus we have:

\begin{prop}\label{2gurka}
Assume that $J\subset M_0$ is a homogeneous submodule 
and \eqref{gcomplex} a free resolution of $M_0/J$ of minimal length, and let
$R$ be  the associated Noetherian residue current.
Let $\phi$ be a holomorphic section of  $E_0\otimes \Ok(r)$ that lies
generically in the image of $f_1\colon E_1\otimes\Ok(r)\to E_0\otimes\Ok(r)$. 
If either

(i)  $N\le n$ 

\noindent  or

(ii) $r\ge \max_j(d^{n+1}_j)-n$,  

\noindent then
$f_1\psi=\phi$ has a global holomorphic solution if (and only if)
$R\phi=0$. 
\end{prop}

Let \eqref{gcomplex} be  any complex and let $R$ be the associated residue current.
If we in addition assume that \eqref{gcomplex} has length at most $n+1$, then
by a similar argument as above it follows that \eqref{gcomplex} is exact if
and only if $R^\ell=0$ for all $\ell\ge 1$, i.e.,
if and only if \eqref{bula} is exact, cf., Theorem~\ref{main1}.

\begin{remark}%%[Some remarks]
The minimal lenght $N$ of a resolution \eqref{gcomplex}
is equal to  $n+1-\depth(M_0/J)$
by the  Auslander-Buchsbaum theorem, see \cite{Eis1}.
The condition (i) is equivalent to that $\depth(M_0/J)\ge 1$ which 
means that  $M_0/J$ contains a nontrivial nonzerodivisor.
If $J$ is defined by a complete intersection,
then the condition (i) is fulfilled. 
 Also if $Z$ is discrete and all the zeros are of
first order, then $\depth S/J=1$, see \cite{Eis2}, so that (i) holds.

\smallskip
The least possible value of $r$  in $(ii)$  is closely related to the
degree of regularity of $J$, see, e.g., \cite{Eis2}.  
An estimate of the regularity for zerodimensional ideals
is given in \cite{Shiff}. See \cite{BayStil} for a general
criterion for a given degree of regularity.
See also Remark~\ref{potta} below.
\end{remark}%%\hspace{3cm} $\square$

%%In analogy with  Proposition~\ref{main1}  the exactness
%%of \eqref{gcomplex} is related to the vanishing of $R$:

\section{Noetherian residue currents for   polynomial  ideals}\label{rotta2}

We will now use the results from the previous section to
obtain Noetherian residue currents for 
(sheaves induced by) polynomial modules in $\C^n$.
Let $z'=(z_1,\ldots,z_n)$ be the standard coordinates in $\C^n$
that we identify  with $\U_0=\{[z]\in\P^n;\ z_0\neq 0\}$,
where $[z]=[z_0,\ldots,z_n]$ are the usual homogeneous coordinates on $\P^n$.
Let $F_1$ be a $\Hom(\C^{r_1}, \C^{r_0})$-valued polynomial in $\C^n$, whose columns 
$F^1,\ldots, F^{r_1}$ have  (at most)  degrees
$d^1_1,\ldots, d^1_{r_1}$ and let $J$ be the submodule of
$\C[z_1,\ldots,z_n]^{r_0}$ generated by $F^1,\ldots, F^{r_1}$.
After the homogenizations  $f^k(z)=z_0^{d_k^1}F^k(z'/z_0)$ 
we get an $r_0\times r_1$-matrix $f_1$ whose columns are $d^1_k$-homogeneous forms in $\C^{n+1}$;
thus a graded mapping 
$$
f_1\colon S(-d_1^1)\oplus\cdots\oplus S(-d_{r_1}^1)\to S^{\oplus r_0}.
$$
%%%%
Extending to a graded resolution (of minimal length)
\eqref{gcomplex}
we obtain  a Noetherian residue current $R$ for the sheaf generated by $f_1$
and an associated current $U$.  In the trivializations
in $\C^n\simeq\U_0$, described  in Example~\ref{badlakan},
the component $R_k$ of $R$  is  the matrix
$(\hat R^0_k)_{ij}$. In the same trivializations
$U_k^\ell$ corresponds to  a matrix $(\hat U^\ell_k)_{ij}$.
Moreover, the mappings $f_k$ correspond to the matrices $F_k$
that are just the dehomogenizations of the matrices $f_k$
in \eqref{plex}.

\smallskip
If $\Phi$ is an  $r_0$-tuple of polynomials in $\C^n$ and there is a
tuple $\Psi$ of  polynomials  such that 
$\Phi=F_1\Psi$  in $\C^n$ then clearly $R\Phi=0$.
Conversely, if $R\Phi=0$ in $\C^n$
(and the equation is locally solvable generically) we know that $\Phi$ is in the
sheaf generated by $F_1$ and hence by Cartan's theorem there is a polynomial solution
to $F_1\Psi=\Phi$. However, we now have a  procedure to find such a $\Psi$:
Take a homogenization $\phi(z)=z_0^r\Phi(z'/z_0)$ for some $r\ge\deg \Phi$.
The condition $R\Phi=0$ in $\C^n$   means that $R\phi=0$ outside the hyperplane at infinity,
so if $r$ is large enough, $R\phi=0$ on $\P^n$.
Now  Proposition~\ref{2gurka} applies if either $r$ is so large that
condition (ii) is fulfilled, or  if  the length of the resolution is less than $n+1$.
If $r$ is chosen large enough we thus have  a holomorphic section $\psi$ of
$E_1\otimes\Ok(r)$ such that $f_1\psi=\phi$. 
After dehomogenization we get the desired polynomial solution $\Psi=(\Psi^j)$
to $F_1\Psi=\sum F^j\Psi^j=\Phi$,  and $\deg F^j\Psi^j\le r$.
%%%
It is well-known that in the worst  case  the final degree has to be  doubly exponential;
at least  $d^{2^{(n/10)}}$,  if  $d$ is the degree of $F_1$, see \cite{MM}.

\begin{remark} \label{potta}
The final degree is essentially depending on the maximal polynomial
degree in the resolution, and it is known to be at worst like $(2d)^{2^n-1}$ if $d$ is the
degree of the generators, see \cite{BayMum}.
\end{remark}

We proceed  with a  result   where we have optimal 
control of the degree of the solution; it is 
a  generalization of Max Noether's classical theorem,  \cite{Noe};
see also  \cite{GH}.

\begin{thm}
Let $F^1,\ldots, F^{r_1}$ be $r_0$-columns of polynomials in $\C^n$
and let $J$ be the homogeneous submodule  of $M_0=S^{\oplus r_0}$ defined by  the homogenized forms
$f^1,\ldots, f^{r_1}$. Furthermore, assume that the quotient module
$M_0/J$  is Cohen-Macaulay and   that no irreducible component of $Z$ is
contained in the hyperplane at infinity. 
If $\Phi$ belongs to the submodule $\tilde J\subset\C[z']^{r_0}$  
generated by $F^1,\ldots, F^{r_1}$, then there are tuples of polynomials
$\Psi^j$ with $\deg (F^j\Psi^j)\le \deg \Phi$ such that
$F^1\Psi^1+\cdots +F^{r_1}\Psi^{r_1}=\Phi$.
\end{thm}

\begin{proof}[Sketch of proof]
We follow the procedure described above. 
Assume that $\codim M_0/J=p$. The Cohen-Macaulay assumption means that
$\dim  M_0/J=\depth M_0/J=n+1-\codim M_0/J$. By the Auslander-Buchsbaum theorem
therefore we can choose a
resolution \eqref{gcomplex}  of $M_0/J$ of length $p$, see \cite{Eis2}.  
Moreover all irreducible components of $Z$ have codimension $p$. 
We choose $r=\deg \Phi$.
Since $\Phi$ is in the ideal in $\C^n$ we have that $R\phi=0$ in $\C^n$.
By Proposition~\ref{huvud},  $R=R_p$ and  
since  $Z$ has no component contained in the hyperplane at infinity, we can copy the argument
in  the proof of Theorem~1.2 in \cite{A4} and conclude that $R\phi=0$ in $\P^n$.
Since $p<n+1$,   cf., Proposition~\ref{2gurka},
we can find a holomorphic section  $\psi$ of $E_1\otimes\Ok(r)$ 
such that $f_1\psi=\phi$. After dehomogenization we get the desired solution $\Psi$.
\end{proof}

We conclude this section with an explicit integral formula that provides  a realization
of the  membership of 
$\Phi$ in  $J\subset \C[z_1,\ldots,z_n]^{r_0}$; for simplicity we assume that
the matrix $F_1=(F^1,\ldots, F^{r_1})$ is generically surjective,
i.e.,  has generic rank $r_0$.
From now on we write  $z$ rather than $z'$.  
It is easy to see that one can choose Hefer matrices  of forms $H^\ell_k$ satisfying
\eqref{Hdef} (with $f_k$ replaced by $F_k$) that are polynomials in both $z$ and $\zeta$;
in fact, the explicit formula in Section~4  in \cite{A7}  when applied to polynomials
will produce polynomials.
%%%
Notice that 
$$
g=\frac{1+\langle\bar\zeta, z\rangle }{1+|\zeta|^2}+\frac{i}{2\pi}\partial\dbar\log(1+|\zeta|^2)
$$
is a weight in $\C^n$ with respect to the point $z$,
cf., Section~\ref{intf}. 
%%Indeed, $g$  is  equal to
%%$
%%1-\nabla_{\zeta-z}
%%\partial\log(1+|\zeta|^2)/2\pi i.
%%$
Since 
$g^{\mu}=\Ok(1/|\zeta|^\mu)$
for fixed $z$  and $H^\ell$ consists of polynomials,  it follows that 
\begin{equation}\label{salsa}
g^{\mu}\w H^0  R, \quad g^{\mu}\w H^1 U
\end{equation}
have current extensions to $\P^n$ if $\mu$ is large enough,
cf., Example~\ref{badlakan}.
Let $\chi_k(\zeta)=\chi(|\zeta|/k)$, where $\chi(t)$ is a cutoff function that
is $1$ for $t<1$ and $0$ for $t>2$. 
If $\mu$ is sufficiently large,
depending on the order at infinity of $R$ and $U$, we have that
\begin{multline}\label{lim}
\chi_k g^\mu\w H^0  R\to g^\mu\w H^0 R,\quad \dbar\chi_k \w g^\mu\w H^0R\to 0, \\
\chi_k g^\mu\w H^1  U\to g^\mu\w  H^1 U,\quad \dbar\chi_k \w g^\mu\w H^1 U\to 0,
\quad 
k\to \infty.
\end{multline}
Let 
$
g_k=\chi_k-\dbar\chi_k\w s/\nabla_{\zeta-z}s,
$
where $s$ is the $(1,0)$-form in \eqref{vikt}.
Then $g_k\w g^{\mu+m}$ is a compactly supported weight with respect to $z$ if $k>|z|$,
and hence   we have  the representation (writing $F$ rather than $F_1$)
$$
\Phi(z)=F(z)\int g_k\w g^{\mu+m}\w H^1 U\Phi+\int g_k\w g^{\mu+m}\w H^0 R\Phi.
$$
Notice   that 
$$
\Big(\frac{1+\langle\bar\zeta, z\rangle}{1+|\zeta|^2}\Big)^m P(\zeta)
$$
is smooth on $\P^n$ for fixed $z$ if $P$ is a polynomial with
$\deg P\le m$. If we let  $k\to \infty$ we therefore obtain

\begin{thm}\label{polyrep}
Let $F$ be a $r_0\times r_1$-matrix of polynomials in $\C^n$ 
with generic rank $r_0$  and let $J$ be the submodule of
$\C[z_1,\ldots,z_n]^{r_0}$ generated by the columns of $F$.
For each given integer $m$, with the notation above and for a large enough $\mu$, 
we have the polynomial decomposition
\begin{equation}\label{polrep}
\Phi(z)=F(z)\int g^{\mu+m}\w H^1 U\Phi+\int  g^{\mu+m}\w H^0 R\Phi
\end{equation}
of $r_0$-columns $\Phi$ of polynomials with  degree at most $m$,
and the last term vanishes as soon as $\Phi\in J$.
\end{thm}

The integrals here are to be interpreted as the action of currents  on
test functions on $\P^n$.
If $\Phi$ belongs to $J$ thus \eqref{polrep} provides a
realization of the membership,
expressed in terms of the current $U$ and the Hefer forms.

%%\smallskip
%%Since  $R$ is Noetherian (for the module sheaf  generated by  $J$)
%%the second integral in \eqref{polrep} vanishes if 
%%$\Phi$ belongs to  $J$;   in that case
%%the first term  thus provides  a  realization of the membership,
%%expressed in terms of the current $U$ and the Hefer forms.

\section{The fundamental principle}\label{fund}

Let $E_1$ and $E_0$ be trivial bundles,  let
 $F$ be a $\Hom(E_1,E_0)$-valued polynomial of generic rank $r_0=\rank E_0$
and let $F^T$ be the transpose of $F$.
Furthermore,   let $K$ be the closure of an open strictly convex bounded  
domain with smooth boundary in $\R^n$ containing the origin.
The fundamental 
principle  of Ehrenpreis and Palamodov states that every homogeneous
solution to the system of equations $F^T(D)\xi=0$,  $D=i\partial/\partial t$,
on $K$ is a superposition of exponential solutions   
with frequencies in the algebraic set $Z=\{z;\ \rank F(z)<r\}$. 
Following the ideas in \cite{BP} we can produce a residue version of the
fundamental principle.

\smallskip

Let $\rho(\eta)$ be the support function 
$
\sup_{t\in K} \langle\eta,t\rangle
$
for $K$ but smoothened  out in a neighborhood of the origin
in $\R^n$. Since $\rho$ is smooth i $\R^n$  and $1$-homogeneous outside a neighborhood
of the origin, all its derivatives are bounded. Let
$$
\rho'(\eta)=( \partial\rho/\partial\eta_1,\ldots, \partial\rho/\partial\eta_n).
$$
We extend to complex arguments
$\zeta=\xi+i\eta$ by letting $\rho(\zeta)=\rho(\eta)$ and $\rho'(\zeta)=\rho'(\eta)$.
Then $\rho'$ maps $\C^n$ onto $K$, see \cite{BP}.
%%Recall that $\nabla_{\zeta-z}=\delta_{\zeta-z}-\dbar$.  Notice that  
%%$$
%%\nabla_{\zeta-z}\partial\rho/2\pi i=-\frac{i}{2}\langle\zeta-z,\rho'(\zeta)\rangle+
%%\partial\dbar\rho(\zeta)/2\pi i.
%%$$
The convexity of $\rho$ implies that
\begin{equation}\label{brass}
e^{\rho(\zeta)}\big|e^{i\langle\rho'(\zeta),\zeta-z\rangle}\big|\le e^{\rho(z)}.
\end{equation}

We  are  to modify  the decomposition \eqref{polrep} to allow
entire functions $h$ with values in $E_0$ satisfying an estimate like
\begin{equation}\label{pally}
|h(z)|\le C(1+|z|)^M e^{\rho(z)}
\end{equation}
for some, from now on,  fixed natural number $M$.
We will use the same notation as in the previous section.
First we introduce a new weight.
%%except that  writing  $F$ rather than $F_1$.
%% but
%%for simplicity we write $R$ and $U$ rather than  $\hat R$ and $\hat U$.
%%and
%%still assume that $R$ is a Noetherian residue current for $J$
%%and $H^\ell_k$ are polynomial-valued Hefer matrices etc.
%%
\begin{lma} The form
$$
g'=e^{i\langle\rho'(\zeta),\zeta-z\rangle+\frac{i}{\pi}\partial\dbar\rho}
=e^{i\langle\rho'(\zeta),\zeta-z\rangle}\sum_{\ell\ge 0}\Big(\frac{i}{\pi}\partial\dbar\rho\Big)^\ell/\ell!
$$
is a weight for each fixed $z\in\C^n$.
\end{lma}

\begin{proof}
Since  $\partial\rho/\partial\zeta_k=-(i/2)\rho'_k(\zeta)$,
$$
\gamma=i\langle\rho'(\zeta), \zeta-z\rangle+\frac{i}{\pi}\partial\dbar\rho(\zeta)=
\nabla_{\zeta-z}\frac{-\partial\rho}{\pi i}
$$
is $\nabla_{\zeta-z}$-closed and $\gamma_{0,0}(z)=0$. 
Thus $\gamma$ and  $e^\gamma$ are  weights.
\end{proof}

It follows from \eqref{brass} that
$$
g^\mu\w g'\w H^1Uh, \quad  g^\mu\w g'\w H^0 Rh
$$
will  vanish to a given finite  order at infinity if $\mu$ is large enough
and $h(\zeta)$ satisfies \eqref{pally}.
Therefore, if $\mu$ is large enough, using the
 compactly supported weights $g_k$  and arguing
 as in the proof of Theorem~\ref{polyrep},  we obtain  the decomposition
\begin{equation}\label{paranot}
h(z)=F(z)\int g'\w g^\mu\w  H^1 Uh +\int g'\w g^\mu \w H^0R h
=F Th +Sh
\end{equation}
for all entire $h$ satisfying \eqref{pally}.
Furthermore,   $Sh$ vanishes if $h=Fq$ for some holomorphic $q$,
and in view of  \eqref{brass}, 
both $Th$ and $Sh$ satisfy \eqref{pally} for some
other large number $M'$ instead of $M$.

\smallskip 

Let $\E'(K)$ be the space of distributions in $\R^n$ with support 
contained in $K$ and let
 $\E^{',M}(K)$ denote the subspace of distributions of order at  most
$M$.  For   $\omega\in\E'(K)$ let $\hat\omega(\zeta)=\omega(e^{-i\langle\zeta,\cdot\rangle})$
be its  Fourier-Laplace  transform.
The Paley-Wiener-Schwartz theorem, see \cite{Hor2} Thm~7.3.1, states that
if $\nu\in\E^{',M}(K)$, then 
\begin{equation}\label{est}
|\hat\nu(\zeta)|\le C(1+|\zeta|)^M e^{\rho(\eta)},
\end{equation}
and conversely:  if $h$ is an entire function that satisfies  such an estimate then 
$h=\hat\nu$ for some  $\nu\in\E'(K)$. 

\smallskip

From \eqref{paranot}, applied to $\hat\nu$ for $\nu\in\E^{',M}(K,E_0)$, 
we therefore get mappings
$$ 
\T\colon \E^{',M}(K,E_0)\to \E'(K,E_1),
\quad\quad \S\colon \E^{',M}(K,E_0)\to \E'(K,E_0),
$$
such that
$$
\nu=F(-D)\T\nu+\S\nu,
$$
and $\S\nu=0$ if $\nu=F(-D)\omega$ for some  $\omega\in\E^{'}(K,E_1)$.
%%%
By duality  we have   mappings
$$
\T^*\colon \E(K, E_1^*)\to C^M(K,E_0^*),\quad
\S^*\colon \E(K, E_0^*)\to C^M(K,E_0^*)
$$
and they satisfy
\begin{equation}\label{daggmask}
\xi =\T^* F^T(D)\xi +\S^*\xi, \quad \xi\in\E(K,E_0^*).
\end{equation}

\begin{thm} Suppose that $M\ge\deg F$.
If $\xi\in\E(K,E_0^*)$, then $\S^*\xi\in C^M(K,E_0^*)$
satisfies $F^T(D) \S^*\xi=0$.
If in addition $F^T(D) \xi=0$, then $S^*\xi=\xi$.
Moreover, we have the explicit formula
\begin{equation}\label{skan}
\S^*\xi(t)=\int_\zeta R^T(\zeta)\alpha^T(\zeta,D)\xi(\rho') e^{-i\langle\zeta,t-\rho'\rangle}
\w  e^{\frac{i}{\pi}\partial\dbar\rho},
\end{equation}
where $\alpha^T(\zeta,D)\xi(\rho')$ is the result when  replacing each
occurrence of $z$ in $\alpha^T(\zeta,z)$ by $D$, letting it act on $\xi(t)$ and evaluating
at the point $\rho'(\zeta)$.
\end{thm}

Thus $\S^*$ is a projection onto the space of homogeneous solutions.

\smallskip
Recall that $\rho'\in K$. Also notice that 
$\Re -i\langle\zeta,t\rangle=\langle \eta,t\rangle\le\rho(\eta)$
if $t\in K$, so combined with \eqref{brass} we get that 
$$
\Re -i\langle\zeta,t-\rho'(\zeta)\rangle\le 0,\quad  t\in K
$$
(for $\zeta$ outside a neighborhood of $0$).
Therefore  the integral in \eqref{skan} has meaning  if $\mu$ is large enough.

\begin{proof}
Suppose that $M\ge\deg F$. Then for $\omega\in\E^{',M-\deg F}(K,E_1)$ we have
\begin{equation}\label{duggmask}
\omega.F^T(D)\S^*\xi=F(-D)\omega.\S^*\xi=\S(F(-D)\omega).\xi=0
\end{equation}
since $\tau=F(-D)\omega\in\E^{',M}(K,E_0)$ so that $\S\tau=0$.
From  \eqref{duggmask}  the first statement now follows.
The second one follows immediately from \eqref{daggmask}. 

It remains to prove \eqref{skan}. The argument is very similar to the proof
of Theorem~2 in \cite{BP} so we only sketch it.
To begin with we have
%%%
\begin{equation}\label{allo}
S\hat\nu(z)=
\int_\zeta \alpha(\zeta,z)R(\zeta)\hat\nu(\zeta) e^{i\langle\zeta-z,\rho'(\zeta)\rangle}
\w e^{\frac{i}{\pi}\partial\dbar\rho}
\end{equation}
where $\alpha(\cdot,z)=g^\mu\w H^0$  is a polynomial in  $z$. 
Let $\delta_t$ be the Dirac measure at $t\in K$. 
Then, letting $T$ denote transpose of matrices, we have
\begin{multline*}
\S^*\xi(t)=\delta_t.\S^*\xi=(\S\delta_t.\xi)^T=\\
\frac{1}{(2\pi)^n}\int_x\int_\zeta
R^T(\zeta)\alpha^T(\zeta,x)e^{-i\langle x,\rho'\rangle}\hat\xi(-x)
e^{-i\langle\zeta,t-\rho'\rangle}\w e^{\frac{i}{\pi}\partial\dbar\rho}.
\end{multline*}
As  in \cite{BP} one can verify that it is legitimate to interchange the order
of integration, and then \eqref{skan} follows by  Fourier's inversion
formula.
\end{proof}

\begin{cor}
For any solution $\xi\in\E(K,E_0^*)$ of $F^T(D)\xi=0$, there are smooth forms
$A_k(\zeta)$ with values in $E_k^*$
such that
\begin{equation}\label{godnatt}
\xi(t)=\int_\zeta \sum_k R^T_k(\zeta) A_k(\zeta) e^{-i\langle\zeta,t-\rho'(\zeta)\rangle}.
\end{equation}
Conversely, for any such smooth forms $A_k(\zeta)$ with sufficient polynomial decay at infinity
the  integral \eqref{godnatt} defines a homogeneous solution.
\end{cor}

The last statement follows just by applying $F^T(D)$ to the integral and using that
$F^T(\zeta)R^T=0$.

\begin{remark}
In case $F$ defines a complete intersection, formulas similar to
\eqref{godnatt} were  obtained in \cite{BP} and \cite{Y}.
In  \cite{BP} is assumed, in addition, that $F^T(D)$ is hypoelliptic; then
one can avoid the polynomial weight factor $g^\mu$ and so the resulting formula is 
even simpler. See also  \cite{BY00} and \cite{BGVY}.
\end{remark}

\begin{ex}[A final example] The ideal $(z_1^2,z_1 z_2)$ corresponds to the
system  
$$
\frac{\partial^2}{\partial t_1^2} \xi(t)=0,
\  \frac{\partial^2}{\partial t_1 \partial t_2} \xi(t)=0.
$$
In view of \eqref{godnatt} and  Example~\ref{enklaste}, the solutions
are precisely the functions that can be written 
\begin{multline*}
\xi(t)=
\int_z 
\Big[\frac{1}{z_2}\Big]\dbar\Big[\frac{1}{z_1}\Big]\wedge A_1(z)~dz_2 \wedge d\bar z_1 \wedge d\bar z_2
~e^{-i(z_1 t_1+ z_2t_2)}+\\
\int_z 
\dbar \Big[\frac{1}{z_1^2}\Big ]\wedge
\dbar \Big[\frac{1}{z_2}\Big ]
\wedge A_2(z) ~d\bar z_1\wedge d\bar z_2
~e^{-i(z_1 t_1+ z_2t_2)},
\end{multline*}
for smooth functions  $A_1$ and $A_2$ with appropriate growth.
It is easily checked  directly to be the general solution, since
the first integral is a quite arbitrary function $C(t_2)$ whereas
the second integral is an arbitrary   polynomial $C_1 +C_2t_1$.
\end{ex}

\def\listing#1#2#3{{\sc #1}:\ {\it #2},\ #3.}

%%Kvarstaende fragor:   

%%Bjorkreferens   fundamentalprincipen?

%%god referens for seminaalytiska mgder

\end{document}